\DocumentMetadata{
    lang        = en-US,
    pdfversion  = 1.7,
    pdfstandard = ua-1,
    tagging     = on,
    tagging-setup = {math/mathml/luamml/load=false}
}
\documentclass[11pt,a4paper]{article}

\usepackage[a4paper,top=2.35cm,bottom=2.6cm,left=3.35cm,right=3.2cm]{geometry}
\usepackage[T1]{fontenc}
\usepackage[american]{babel}

\usepackage[sfdefault,light]{FiraSans}
\usepackage{newtxsf}

\usepackage{xcolor}
\usepackage{microtype}
\usepackage{setspace}
\usepackage{tikz}
\usepackage{eso-pic}
\usepackage{array}
\usepackage{hyperref}
\usepackage{tabularx}

\definecolor{Ink}{HTML}{20242A}
\definecolor{Slate}{HTML}{606B75}
\definecolor{Mist}{HTML}{F4F7F8}
\definecolor{MistBlue}{HTML}{EAF2F6}
\definecolor{LeidenNavy}{HTML}{0B243B}
\definecolor{LeidenBlue}{HTML}{174A6E}
\definecolor{LeidenTeal}{HTML}{2A7E80}
\definecolor{LeidenGold}{HTML}{BD9142}
\definecolor{WarmGrey}{HTML}{E7E2D8}
\definecolor{RuleGrey}{HTML}{D9DEE3}

\color{Ink}

\hypersetup{
    colorlinks=true,
    linkcolor=LeidenBlue,
    urlcolor=LeidenTeal,
    citecolor=LeidenBlue,
    pdftitle={Leiden Declaration on Artificial Intelligence and Mathematics},
    pdfauthor={Leiden Declaration working group},
    pdfsubject={Recommendations on artificial intelligence and mathematics research},
    pdfkeywords={artificial intelligence, mathematics, research ethics, open science},
    pdfcreator={LaTeX},
    pdfdisplaydoctitle=true,
    pageanchor=false,
    bookmarksopen=true,
    bookmarksnumbered=false
}

\newcommand{\weblink}[2]{\href{#1}{\textcolor{LeidenTeal}{#2}}}
\newcommand{\UNESCOOpenScience}{\weblink{https://www.unesco-floods.eu/wp-content/uploads/2022/04/379949eng.pdf}{UNESCO Recommendation on Open Science}}

\makeatletter
\renewcommand\section{\@startsection{section}{1}{0pt}%
  {2.1em}%
  {0.5em}%
  {\sffamily\Large\bfseries\color{LeidenNavy}}}
\renewcommand\subsection{\@startsection{subsection}{2}{0pt}%
  {0.45em}%
  {0.15em}%
  {\sffamily\large\bfseries\color{LeidenBlue}}}
\newcommand{\SectionDecorRule}{%
  \vspace{-1.35em}%
  \tagmcbegin{artifact=layout}%
  \noindent\textcolor{LeidenGold}{\rule{\linewidth}{1pt}}\par
  \tagmcend
  \vspace{-0.32em}%
}
\let\WCAGBaseSection\section
\renewcommand{\section}{%
  \@ifstar{\WCAGSectionStar}{\@ifnextchar[{\WCAGSectionOpt}{\WCAGSectionNoOpt}}%
}
\newcommand{\WCAGSectionStar}[1]{\WCAGBaseSection*{#1}\SectionDecorRule}
\newcommand{\WCAGSectionOpt}[2][]{\WCAGBaseSection[#1]{#2}\SectionDecorRule}
\newcommand{\WCAGSectionNoOpt}[1]{\WCAGBaseSection{#1}\SectionDecorRule}
\makeatother

\newcommand{\CoverBackground}{%
\begin{tikzpicture}[artifact,remember picture,overlay]

    \fill[MistBlue]
        (current page.north west)
        rectangle ([yshift=-9.7cm]current page.north east);

    \fill[LeidenNavy]
        (current page.north west)
        rectangle ([xshift=0.36cm]current page.south west);

    \fill[LeidenGold]
        ([xshift=0.36cm]current page.north west)
        rectangle ([xshift=0.48cm]current page.south west);

    \foreach \x/\y/\r in {
        -1.2/-1.4/3.9,
        -0.4/-5.3/2.5,
        -4.8/-2.2/1.4,
        -6.5/-6.6/0.85
    }{
        \draw[LeidenTeal,opacity=0.16,line width=0.9pt]
            ([xshift=\x cm,yshift=\y cm]current page.north east) circle (\r cm);
    }

    \foreach \x/\y in {
        -3.1/-2.8,
        -4.4/-4.1,
        -2.2/-5.1,
        -6.0/-3.0,
        -5.5/-6.2
    }{
        \fill[LeidenGold,opacity=0.42]
            ([xshift=\x cm,yshift=\y cm]current page.north east) circle (2pt);
    }

    \draw[LeidenTeal,opacity=0.32,line width=0.6pt]
        ([xshift=-3.1cm,yshift=-2.8cm]current page.north east) --
        ([xshift=-4.4cm,yshift=-4.1cm]current page.north east) --
        ([xshift=-2.2cm,yshift=-5.1cm]current page.north east) --
        ([xshift=-5.5cm,yshift=-6.2cm]current page.north east) --
        ([xshift=-6.0cm,yshift=-3.0cm]current page.north east) --
        ([xshift=-3.1cm,yshift=-2.8cm]current page.north east);

\end{tikzpicture}
}

\newcommand{\PageAccent}{%
\begin{tikzpicture}[artifact,remember picture,overlay]
    \fill[LeidenNavy]
        (current page.north west)
        rectangle ([xshift=0.18cm]current page.south west);
    \fill[LeidenGold]
        ([xshift=0.18cm]current page.north west)
        rectangle ([xshift=0.25cm]current page.south west);
    \node[anchor=base west,inner sep=0pt] at ([xshift=3.343cm,yshift=-1.345cm]current page.north west)
        {\sffamily\small\textcolor{Slate}{Leiden Declaration}};
    \node[anchor=base east,inner sep=0pt] at ([xshift=-3.193cm,yshift=-1.345cm]current page.north east)
        {\sffamily\small\textcolor{Slate}{Artificial Intelligence and Mathematics}};
    \node[anchor=center] at ([xshift=0.075cm,yshift=1.67cm]current page.south)
        {\sffamily\small\textcolor{Slate}{\thepage\ /\ 11}};
\end{tikzpicture}
}

\newcommand{\standfirst}[1]{%
\vspace{0.2em}
\noindent
\raisebox{10.79em}[0pt][0pt]{\makebox[0pt][l]{%
\begin{tikzpicture}[artifact,baseline=(standfirstbox.north)]
\node[
    fill=Mist,
    draw=RuleGrey,
    rounded corners=8pt,
    inner sep=16pt,
    minimum width=\dimexpr0.88\textwidth+32pt\relax,
    minimum height=20.68em,
    anchor=north west
] (standfirstbox) at (0,0) {};
\end{tikzpicture}%
}}%
\begin{minipage}{\dimexpr0.88\textwidth+32pt\relax}
\vspace{0.925em}
\hspace{16pt}%
\begin{minipage}{\dimexpr\linewidth-32pt\relax}
\sffamily\large\setstretch{1.085}\color{LeidenNavy}%
#1\par
\end{minipage}
\vspace{1.175em}
\end{minipage}
\par
\vspace{0.55em}
}

\newcounter{principlebookmark}
\newcommand{\principle}[2]{%
\par\vspace{1.0em}
\refstepcounter{principlebookmark}%
\pdfbookmark[2]{#1. #2}{principle.\theprinciplebookmark}%
\noindent
\raisebox{2.1em}[0pt][0pt]{\makebox[0pt][l]{%
\begin{tikzpicture}[artifact,baseline=(principlebox.north)]
\node[
    fill=white,
    draw=RuleGrey,
    rounded corners=8pt,
    inner sep=0pt,
    minimum width=0.88\textwidth,
    minimum height=3.65em,
    anchor=north west
] (principlebox) at (0,0) {};
\end{tikzpicture}%
}}%
\begin{minipage}{0.7744\textwidth}
\vspace{1.0em}
\begin{minipage}[t]{0.12\linewidth}
    \centering
    \sffamily\bfseries\Large\textcolor{LeidenGold}{#1}
\end{minipage}
\begin{minipage}[t]{0.85\linewidth}
    \sffamily\bfseries\large\textcolor{LeidenNavy}{#2}
\end{minipage}
\vspace{0.95em}
\end{minipage}
\par\vspace{0.55em}
}

\newcommand{\numberedpoint}[2]{%
\par\noindent\hangindent=1.15em\hangafter=1
\makebox[1.15em][l]{#1.}#2\par\vspace{-0.18em}
}

\begin{document}


\thispagestyle{empty}
\AddToShipoutPicture*{\CoverBackground}

\vspace*{2.65cm}

\vspace{0.85cm}

{\sffamily\bfseries\fontsize{36}{42}\selectfont\color{LeidenNavy}
Leiden Declaration\\
on Artificial Intelligence\\
and Mathematics
\par}

\vspace{1.9cm}


\vspace{1.15cm}

{\color{LeidenGold}\rule{0.28\textwidth}{2.2pt}}

\vfill

\begin{minipage}{1.0\textwidth}
\sffamily\Large\color{Slate}
This declaration calls for action to address the challenges posed by the use of artificial intelligence within mathematics research.
\end{minipage}

\vspace{2.0cm}

\begin{flushright}
\sffamily
\textcolor{Slate}{June 2026}\\
\vspace{1.15cm}
\weblink{https://doi.org/10.5281/zenodo.20302944}{DOI: 10.5281/zenodo.20302944}
\end{flushright}

\newpage

\AddToShipoutPicture{\PageAccent}


\section{Preamble}

Technological developments have repeatedly transformed the practice of mathematics. Recent artificial intelligence technologies, including symbolic and neural methods for the generation and formalization of mathematics, may already have initiated a significant chapter in this long history. Among researchers, artificial intelligence has produced a wide range of reactions: enthusiasm for its potential to yield new discoveries; intimidation by the pace of developments; indifference to these rapid changes; and concern for the implications, both for mathematics and in wider society.

Mathematicians have a choice about whether and how to adopt artificial intelligence in the conduct of their research. They also have a responsibility to ensure the continued flourishing of the discipline. This Declaration calls upon mathematicians to exercise this responsibility, and provides recommendations for individuals, institutions, government, and
industry.

Although we adopt the perspective of mathematical research, much of what we write applies equally to other aspects of mathematics.   This includes work in the broader mathematical sciences, education, mentoring, publishing, funding, science policy, and use of mathematics in the wider world.

The Declaration is conceived in solidarity with other research endeavors and creative professions facing similar challenges, both within and beyond academia. It complements other calls  for action such as the \weblink{https://phsj.org/wp-content/uploads/2007/10/Uppsala-Code-of-Ethics-for-Scientists.pdf}{Uppsala Code of Ethics for Scientists}, the \weblink{https://sfdora.org/read/}{San Francisco Declaration on Research Assessment}, the \weblink{https://www.unesco-floods.eu/wp-content/uploads/2022/04/379949eng.pdf}{UNESCO Recommendation on Open Science}, and the \weblink{https://www.gov.uk/government/publications/universal-ethical-code-for-scientists}{UK Universal Ethical Code for Scientists}. The \weblink{https://www.mathunion.org/cop}{International Mathematical Union Committee on Publishing,} the \weblink{https://www.siam.org/media/b03hwuwe/siam-report-ai-task-force.pdf}{Society for Industrial and Applied Mathematics}, and the \weblink{https://www.ams.org/about-us/ai-summary}{American Mathematical Society} have  also produced related material.\looseness-1

\subsection{About our values}

We base our recommendations on what we take to be characteristic values of mathematical research that we have a joint interest in preserving. Among these are the following:

\numberedpoint{1}{There are many reasons to pursue mathematical research, ranging from intellectual curiosity to a desire to solve practical and societal problems. Underlying much of mathematics is the activity of proof. Mathematical proofs are regarded as conferring the highest degree of certainty to their conclusions, as well as imparting understanding of why their conclusions are true. These characteristics of proof support the scientific integrity  of mathematics.}

\numberedpoint{2}{Results are attributable to specific authors who take credit for their discovery and assume responsibility for their correctness. These principles ground the merit-based standards to which we aspire in mathematical research.}

\numberedpoint{3}{Mathematical arguments are regarded as transparent and subject to independent verification. They may be extremely long or difficult, but in principle no proprietary knowledge or equipment should be required to understand them.}

\numberedpoint{4}{Mathematicians share a concern for proper evaluation of mathematical work relative to shared standards of depth, difficulty, and significance.}

\numberedpoint{5}{Mathematics produces not only a body of results, but also understanding, clarity, and judgment among the communities of mathematicians who have shaped them, often in the context of their own autonomously guided research.  This expert knowledge is essential, both to effectively use mathematics, and to continue to articulate new and significant research questions. A key source of strength of the discipline has long been the autonomous shaping of the direction of research and the methods used to pursue it.}

These characteristics of mathematics as a subject matter are also compatible with understanding mathematics as a human practice, and its place in the world. As mathematicians, and also as inhabitants of a shared world, we have a duty to care for other people and our environment.

\subsection{Potential threats}

Recent developments in  artificial intelligence threaten each of these values, often in ways that disproportionately affect students and early-career mathematicians, and hence the long term future of the discipline.

\numberedpoint{1}{Current automated techniques can produce plausible but unreliable (or even incorrect) arguments which are difficult to distinguish from correct mathematical proofs. This applies not only to informal arguments, but also to formalizations, where the difficulty lies in the translation between computer-encoded and human presentations of concepts. These fast-moving developments put our present system of review under increasing pressure, jeopardizing our ability to implement traditional standards for the correctness, transparency, and independent verifiability of proof.}

\numberedpoint{2}{Technologies that draw extensively on the published mathematical commons undermine the traditional system of attribution. Models trained on published works frequently return outputs that do not properly cite the human works they synthesize.
Many current models are also built on data obtained by systematically exploiting licenses and access arrangements that were not made with artificial intelligence in mind, or indeed by simply violating copyright protections.}

\numberedpoint{3}{Technologies which affect the way in which mathematics is practiced may disturb the current system of incentives.  The use of artificial intelligence -- and thus also the sort of problems which it can  address -- may become incentivized for its own sake, disrupting our mechanisms for hiring, funding, and recognition. This disadvantages researchers who do not have access to the technologies or decision-making related to them, or who are unwilling to use technologies controlled by organizations whose values they do not share.}

\numberedpoint{4}{Proper evaluation is endangered if results are communicated through informal channels such as press releases or blog posts, often without any research paper or other disclosure of information necessary for scientific evaluation. This practice seeks publicity for new results on market timelines before the accepted processes of community evaluation in mathematics can take place. In many cases this leads to simplifications in reporting, such as overemphasizing the significance of automated tools and undervaluing the prior human contributions which have made those tools possible. Such oversimplification risks influencing public opinion in a way that not only damages perceptions of mathematics, but also misleadingly uses specific mathematical tasks as metrics for the general reasoning capacities of commercial products.\looseness-1}

\numberedpoint{5}{These developments put the autonomy of mathematics under threat.  The increasing involvement of technology companies in mathematical research raises the risk that research questions may come to be prioritized because of their amenability to automated mathematics, rather than expert judgment of their deeper significance.  Indeed, broader understanding of the field may be permanently lost in the process of automation. With university budgets under pressure, this reshaping also changes professional incentives in a manner which encourages the collaboration of researchers with technology companies on asymmetric terms. If left unchecked, these trends go beyond threatening researchers' autonomy, affecting the scope and depth of mathematical research itself.}

All of these challenges arise at a moment when the consequences of large-scale investment in artificial intelligence are being widely discussed in regard to warfare, mass surveillance, political disruption, and environmental damage. These raise grave ethical concerns. By failing to act, we run the risk of becoming complicit in the support of technologies which threaten much more than the practice of mathematics.

We thus feel that there is an urgent need for a considered response from the mathematical community. The following constitute brief descriptions of actionable recommendations. We encourage professional organizations to endorse this Declaration, and to add provisions according to their own values, priorities, and governance.


\section{Recommendations for individual mathematicians}

\principle{01}{Disclose tool use} \label{sec:authorship}

Transparently disclose the use of automated tools, including large language models, machine learning systems, proof assistants, and other mathematical software. Include a ``Tool and computational resource disclosure'' section in your papers; many journals, publishers, and professional organizations have already developed guidelines for this, and though the precise form of such a section will necessarily evolve, we encourage authors to live up to the spirit reflected in the \UNESCOOpenScience\ and the \weblink{https://www.nature.com/articles/sdata201618}{FAIR principles}. When acting as a reviewer, abide by publisher guidelines. If the use of artificial intelligence is allowed, be transparent about how you used it, and take responsibility for any significant recommendations you make.

\newpage
\principle{02}{Support the needs of reviewing}

The use of artificial intelligence in preparing papers can introduce material that makes reviewing more demanding. Make it easier for your peers to review your work by disclosing tool use, giving precise and complete references to previous results, and providing formal proofs where feasible and appropriate.

\principle{03}{Adhere to principles of open science}

The international open science movement aims to make scientific research transparent and accessible to all. As mathematical research becomes more reliant on data and software, adhere to principles of open science. See also the \UNESCOOpenScience.

\principle{04}{Retain the responsibility for correctness}

When automated techniques are employed in published mathematical research, the responsibility for the correctness and adequacy of the arguments and results, as well as for the completeness and accuracy of citations to relevant prior work, remains exclusively with the human authors.

\principle{05}{Affirm the humanity of authorship}

Credit and responsibility continue to belong to humans within the mathematical community and should not be given to automated systems. Artificial intelligence may obscure, but does not replace, the collective human labor behind a result.

\principle{06}{Put effort into proper attribution}

The known limitations of automated tools in properly attributing ideas create a corresponding obligation for proactive effort to find and credit the sources  that made a new result possible. Where a satisfactory attribution is not possible, state this explicitly in the publication.

\newpage

\principle{07}{Participate in public discourse}

Mathematicians have a responsibility to support serious science journalism and to engage in public discourse to explain and contextualize artificial intelligence-assisted methods and results. This is particularly important for work within our own subfields, where specialized knowledge is required to assess claims about the depth, difficulty, and significance of results. Moreover, we encourage mathematicians to seek opportunities to cooperate with and support  other researchers and creative professionals facing similar challenges.

\principle{08}{Stay informed about the emerging technologies}

As appropriate to your interests and research, stay informed about the capability of computer-aided mathematical tools. Such understanding is important for informing how our discipline adapts to new technologies and for participating in governance and public discourse.

\principle{09}{Welcome new contributors}

The growing intersection of artificial intelligence and mathematics continues to attract researchers from other disciplines. We welcome this broadening of our community and the range of skills and perspectives these contributors bring. We encourage the mathematical community to actively engage with the broader community, to make our standards and practices explicit and accessible, and to create pathways for meaningful participation. In turn, we ask those entering our field to approach it with respect for our values, while also helping us to adapt and develop them.

\principle{10}{Consider carefully which tools to use}

Some automated tools and their developers  will align with the provisions of this Declaration, while others will not. Consider this when deciding which tools to use, or whether to use them at all. Also consider whether non-proprietary, energy-efficient, or small-scale systems suffice for your task. If not, consider how preservation of the values articulated in this Declaration may be worth a delay in obtaining results.

\newpage

\principle{11}{Evaluate the ethical consequences of your work, and take action accordingly}

Mathematics has led to technology which greatly improves everyday life for many people, yet it also has applications in the development of technology for use in warfare, oppression, mass surveillance, and the undermining of democracy. Evaluate the ethical consequences of your research to the best of your abilities, and if necessary withdraw from harmful work. Only enter into external partnerships which respect the values articulated in this Declaration.

\section{Recommendations for mathematical organizations and \texorpdfstring{\\}{ }not-for-profit research funders}

\principle{01}{Build expertise and plan strategically}

Professional organizations should keep abreast of technical developments and be proactive in making informed recommendations to members and to the broader community. They should work together to guide the development of policy within academic publishing, funding bodies, and government. They should also actively prepare to become involved if major mathematical results are claimed using unconventional means.\looseness-1

\principle{02}{Take the lead on policies for publishing and reviewing}

Professional organizations within mathematics should take a leadership role in developing guidelines in regard to the use of automated techniques in publication and in reviewing. These would include, for example, tool and computational resource disclosure, attribution, rules pertaining to authorship, and codes of conduct consistent with the values of mathematics. These would  supplement and support guidelines already being developed by publishers and journals.

\principle{03}{Maintain standards of rigor}

When establishing policies, demand that results obtained  by automated techniques be held to standards that address the risks raised by those techniques. These might include requiring human descriptions of central arguments obtained by automated tools,  insisting on formal verification when appropriate, cross-checking theoretical and computational results, or external pre-submission review.

\principle{04}{Protect the rights of authors}

Automated mathematics presents new challenges to the rights of authors, and societies should be proactive in the development of sample licensing agreements to protect these rights. In particular, material  should not be used as training data without 
consent, and publishing agreements should allow authors  to opt-out of the use of their work in this way.

\principle{05}{Insist on appropriate publication outlets}

Demand that mathematical results continue to be published in peer-reviewed venues such as journals, proceedings, and books. Informal mechanisms such as press releases or blog posts can provide a valuable supporting role, but they cannot replace peer-review or community scrutiny.

\principle{06}{Support public research laboratories}

Support the formation of university-based, national, or international research laboratories devoted to studying automated mathematics which are administratively and financially independent from industry. Support the use of less resource-intensive technologies accessible to individual researchers.

\principle{07}{Provide frameworks for collaboration}

Mathematicians and academic organizations collaborating with industry often face asymmetries in their bargaining positions, as well as in access to professional support such as legal resources, or advice on intellectual property.   Support researchers in such collaborations by providing access to legal representation, and by facilitating the development of codes of professional practice.

\principle{08}{Align funding with values}

Alignment with the values of this declaration should be taken into account in the evaluation and funding of projects which involve collaboration between academics and industrial partners.

\pagebreak
\section{Recommendations for policymakers in government and \texorpdfstring{\\}{ }elsewhere}

\principle{01}{Protect the rights of authors}

Strengthen legal protections for authors, in line with this declaration.

\principle{02}{Don't believe the hype}

There is currently a strong commercial incentive on the part of the technology industry to overstate the capabilities of their products. Consult with experts, including mathematicians, in forming policy decisions rather than relying on press releases or popular reporting of mathematical results.

\principle{03}{Regulate the artificial intelligence industry}

Recent developments continue to highlight the strong public interest in regulating the technology industry, for example in regard to involvement in military and mass surveillance programs, development of technologies which promote misinformation and undermine democracy, and environmental costs. We stand with others in calling for significantly increased public oversight.

\principle{04}{Invest in public computational infrastructure}

Current events illustrate the need for public alternatives to proprietary technologies, from basic services for online collaboration, to computer clusters for mathematical modeling and machine learning applications. We support the funding of public infrastructure at university, national, and international levels.

\section{Recommendations for commercial artificial intelligence}

While the mathematical community has recognized standing in academic and public policymaking, it has no comparable role in the corporate decision-making that is playing an increasing  role in our discipline.   Nonetheless, recent developments have drawn mathematical work into industrial artificial intelligence efforts in multiple ways. One is through the use of mathematics to advertise the capabilities of commercial artificial intelligence systems in public communications and public relations campaigns. Another is that artificial intelligence developers have increasingly used mathematical publications and formal mathematical libraries as sources of training data -- not only for specialized models for mathematics, but for more general-purpose artificial intelligence.\looseness-1

What currently makes mathematics attractive for general-purpose artificial intelligence development is that the correctness of formalized proofs can be checked automatically, without the need for human oversight. This makes it possible to generate and check vast numbers of problems, both human-authored and computer-generated, to produce an effectively unlimited source of feedback for training artificial intelligence models. The rationale for this strategy often rests on a further assumption: that capabilities developed through mathematical theorem proving will extend to broader general reasoning. Some of the resulting general-purpose models are being commercialized for applications that raise grave ethical concerns, including those named earlier: warfare, oppression, mass surveillance, and the undermining of democracy.

We recognize that industry has offered lucrative jobs, monetary rewards, computing resources, and intellectually stimulating opportunities that some mathematicians have found attractive. This has taken place in an era of underfunding of higher education and precarious academic employment. We also recognize that many mathematicians did not expect their work to become entangled with social and ethical implications of such magnitude, nor to be incorporated into systems used for purposes they may find deeply troubling.

We call on collaborations between mathematicians and industry to abide, at minimum, by the standards we expect of our colleagues and that are described throughout this Declaration. Such collaborations must respect the freedom of conscience of employees or contributors to speak openly about corporate policies and priorities.

\vspace{-0.77em}
\section{About this document}
\vspace{0.77em}
\standfirst{%
In September 2025 the Lorentz Center at Leiden University in the Netherlands hosted a conference entitled \weblink{https://www.lorentzcenter.nl/mechanization-and-mathematical-research.html}{Mechanization and Mathematical Research}. The around 60 participants from 10 countries comprised mathematicians, computer scientists, philosophers, historians, and social scientists, including those with experience in industry and in government. During the eight months following the conference, a smaller working group developed this Declaration, with extensive feedback from the mathematical community.
The Declaration reflects artificial intelligence technologies and mathematical practice as of May 2026. The Declaration can be signed, with authentication through an ORCID identifier, at  \href{https://leidendeclaration.ai}{\mbox{\textcolor{LeidenTeal}{\texttt{https://leidendeclaration.ai}}}}. The working group was convened by Jim Portegies {\tt <j.w.portegies@tue.nl>} who can be contacted for any further information.}

\newpage

\section{Members of the working group}

\vspace{0.5em}
\tagpdfsetup{table/header-rows={1}}
\renewcommand{\arraystretch}{1.1}
\begin{tabularx}{0.88\textwidth}{>{\rule[-1.5\baselineskip]{0pt}{1\baselineskip}}X X}

\sffamily\small Jarod Alper & \sffamily\small University of Washington \\

\sffamily\small Michael Barany & \sffamily\small University of Edinburgh \\

\sffamily\small Alain Chavarri Villarello & \sffamily\small Vrije Universiteit Amsterdam \\

\sffamily\small Sander Dahmen & \sffamily\small Vrije Universiteit Amsterdam \\

\sffamily\small Walter Dean & \sffamily\small University of Warwick \\

\sffamily\small Karthik Ganapathy & \sffamily\small University of California, San Diego \\

\sffamily\small Michael Harris & \sffamily\small Columbia University \\

\sffamily\small David Holmes & \sffamily\small Leiden University \\

\sffamily\small Mateja Jamnik & \sffamily\small University of Cambridge \\

\sffamily\small Steven Kelk & \sffamily\small Maastricht University \\

\sffamily\small Bryna Kra & \sffamily\small Northwestern University \\

\sffamily\small Ursula Martin & \sffamily\small University of Oxford \\

\sffamily\small Bartosz Naskręcki  & \parbox[t]{\linewidth}{\sffamily\small Adam Mickiewicz University\\\sffamily\small Warsaw University of Technology} \\

\sffamily\small Rodrigo Ochigame & \sffamily\small Leiden University \\

\sffamily\small Jim Portegies & \sffamily\small Eindhoven University of Technology\\

\sffamily\small Johannes Schmitt & \sffamily\small ETH Zurich \\

\end{tabularx}

\vfill

\tagmcbegin{artifact=layout}%
\par\noindent\makebox[\textwidth][c]{\textcolor{RuleGrey}{\rule{0.62\textwidth}{0.6pt}}}\par
\vspace{0.8em}
\par\noindent\makebox[\textwidth][c]{%
  \sffamily\small\textcolor{Slate}{Leiden Declaration on Artificial Intelligence and Mathematics}%
}\par
\tagmcend
\end{document}